\documentclass{amsart}

\newcommand{\F}{\mathbb{F}}
\newcommand{\BP}{\mathbb{P}}

\def\calO{\mathcal{O}}

\def\frakp{\mathfrak{p}}
\def\frakq{\mathfrak{q}}

\DeclareMathOperator{\Gal}{Gal}
\DeclareMathOperator{\Jac}{Jac}

\newtheorem{thm} {Theorem}
\newtheorem{lem} [thm]{Lemma}

\newtheorem*{heur}{Heuristic}

\newtheorem*{problem}{Problem}

\begin{document}

\title[Pointless curves]
      {Pointless curves of genus three and four}

\author{Everett W.~Howe}
\address{Center for Communications Research, 
         4320 Westerra Court, 
         San Diego, CA 92121-1967, USA.} 
\email{however@alumni.caltech.edu}
\urladdr{http://www.alumni.caltech.edu/\~{}however/}

\author{Kristin E. Lauter}
\address{Microsoft Research,
         One Microsoft Way,
         Redmond, WA 98052, USA.}
\email{klauter@microsoft.com}

\author{Jaap Top}
\address{Department of Mathematics, 
         University of Groningen, 
         P.O. Box 800, 
         9700 AV Groningen, 
         The Netherlands.}
\email{top@math.rug.nl}

\date{1 March 2004}



\begin{abstract}
A curve over a field $k$ is {\em pointless\/} if it has no $k$-rational
points.  We show that there exist pointless genus-$3$ hyperelliptic curves over
a finite field $\F_q$ if and only if $q\le 25$,
that there exist pointless smooth plane quartics over $\F_q$ if and only if
either $q\le 23$ or $q=29$ or $q=32$, and 
that there exist pointless genus-$4$ curves 
over $\F_q$ if and only if $q\le 49$.
\end{abstract}

\maketitle

\section{Introduction}
\label{S-intro}

How many points can there be on a curve of genus $g$ over a finite
field~$\F_q$?  Researchers have been studying variants of this question
for several decades.
As van der Geer and van der Vlugt write in the 
introduction to their biannually-updated survey of results
related to certain aspects of this subject,
the attention paid to this question is
\begin{quote}
motivated partly by possible applications in coding
theory and cryptography, but just as well by the fact that 
the question represents an attractive mathematical challenge.
\cite{GeerVlugt}
\end{quote}
The complementary question --- how {\em few\/} points can there be
on a curve of genus $g$ over $\F_q$? --- seems to have
sparked little interest among researchers,
perhaps because of the apparent {\em lack\/}
of possible applications for such curves in coding
theory or cryptography.  
But despite the paucity of applications, there 
are still mathematical challenges associated with such curves.
In this paper, we address one of them:
\begin{problem}
Given an integer $g\ge 0$, determine the finite fields $k$ such 
that there exists a curve of genus $g$ over $k$ having
no rational points.
\end{problem}
We will call a curve over a field $k$ {\em pointless\/} if
it has no $k$-rational points.  Thus the problem we propose
is to determine, for a given genus $g$, the finite fields $k$
for which there is a pointless curve of genus~$g$.

The solutions to this problem for $g\le 2$ are known.
There are no pointless curves of genus~$0$ over any finite field;
this follows from Wedderburn's theorem, as is shown 
by~\cite[\S~III.1.4, exer.~3]{Serre:CG}.
The Weil bound for curves of genus $1$ over a finite field,
proven by Hasse~\cite{Hasse}, shows that there are no pointless
curves of genus~$1$ over any finite field.
If there is a pointless curve of genus $2$ over a finite field~$\F_q$ 
then the Weil bound shows that $q\le 13$, and in 1972 Stark~\cite{Stark}
showed that in fact $q < 13$.  For each $q<13$ there do
exist pointless genus-$2$ curves over $\F_q$; a complete list 
of these curves is given in~\cite[Table~4]{MaisnerNart}.

In this paper we provide solutions for the cases $g=3$ and $g=4$.

\begin{thm}
\label{T-genus3}
There exists a pointless genus-$3$ curve over $\F_q$ if and
only if either $q\le 25$ or $q=29$ or $q=32$.
\end{thm}

\begin{thm}
\label{T-genus4}
There exists a pointless genus-$4$ curve over $\F_q$ if and 
only if $q \le 49$.
\end{thm}

In fact, for genus-$3$ curves we prove a statement slightly
stronger than Theorem~\ref{T-genus3}:

\begin{thm}
\label{T-genus3specific}
\ 
\begin{enumerate}
\item There exists a pointless genus-$3$ hyperelliptic curve 
over $\F_q$ if and only if~$q\le 25$.
\item There exists a pointless smooth plane quartic curve over $\F_q$
if and only if either $q\le 23$ or $q=29$ or $q=32$.
\end{enumerate}
\end{thm}

The idea of the proofs of these theorems is simple.  For any 
given genus $g$, and in particular for $g=3$ and $g=4$,
the Weil bound can be used to provide an upper bound for the set
of prime powers $q$ such that there exist pointless curves of 
genus $g$ over $\F_q$.  For each $q$ less than or equal to this bound, we
either provide a pointless curve of genus $g$ or use the techniques
of~\cite{HoweLauter} to prove that none exists.

We wrote above that the question of how few points there can
be on a genus-$g$ curve over $\F_q$ seems to have attracted
little attention, and this is certainly the impression one gets
from searching the literature for references to such curves.
On the other hand, the question has undoubtedly occurred to researchers before.
Indeed, the third author was asked this very question for
the special case $g=3$ by both N.~D.~Elkies and J.-P.~Serre
after the appearance of his joint work~\cite{AuerTop} with Auer.
Also, while it is true that there seem to be no applications for
pointless curves, it {\em can\/} be useful to know whether or
not they exist.  For example, Leep and Yeomans
were concerned with the existence of pointless
plane quartics in their work~\cite{LeepYeomans}
on explicit versions of special cases of the Ax-Kochen theorem.
Finally, we note that Clark and Elkies have recently proven
that for every fixed prime $p$ there is a constant $A_p$ such 
that for every integer $n>0$ there is a curve over $\F_p$
of genus at most $A_p n p^n$ that has no places of degree $n$ or less.

In Section~\ref{S-heuristics} we give the heuristic that 
guided us in our search for pointless curves.  In Section~\ref{S-proof}
we give the arguments that show that there are no pointless curves
of genus $3$ over $\F_{27}$ or $\F_{31}$, 
no pointless smooth plane quartics over $\F_{25}$,
no pointless genus-$3$ hyperelliptic curves over $\F_{29}$ or $\F_{32}$,
and no pointless curves of genus $4$ over $\F_{53}$ or~$\F_{59}$. 
Finally, in Sections~\ref{S-examples3} and~\ref{S-examples4}
we give examples of pointless curves of genus
$3$ and $4$ over every finite field for which such curves exist.

\subsubsection*{Conventions}
By a {\em curve\/} over a field $k$ we mean a smooth, projective,
geometrically irreducible $1$-dimensional variety over~$k$.
When we define a curve by a set of equations, we mean the normalization
of the projective closure of the variety defined by the equations.

\subsubsection*{Acknowledgments}
The first author spoke about the work~\cite{HoweLauter}
at AGCT-9, and he thanks the organizers
Yves Aubry, Gilles Lachaud, and Michael Tsfasman
for inviting him to Luminy and for organizing
such a pleasant and interesting
conference.  The first two authors thank the
editors for soliciting this paper, which made them think
about other applications of the techniques developed 
in~\cite{HoweLauter}.

In the course of doing the work described in this paper
we used the computer algebra system Magma~\cite{magma}.
Several of our Magma programs are available on the web:  start at

\smallskip
\noindent
{\tt http://www.alumni.caltech.edu/\~{}however/biblio.html}
\smallskip

\noindent
and follow the links related to this paper. 
One of our proofs depends on an explicit description of the
isomorphism classes of unimodular quaternary Hermitian forms 
over the quadratic ring of discriminant $-11$.
The web site mentioned above also contains a copy of a text file 
that gives a list of the six isomorphism classes of 
such forms; we obtained this file from the web site

\smallskip
\noindent
{\tt http://www.math.uni-sb.de/\~{}ag-schulze/Hermitian-lattices/}
\smallskip

\noindent
maintained by Rainer Schulze-Pillot-Ziemen.


\section{Heuristics for constructing pointless curves}
\label{S-heuristics}

To determine the correct statements of Theorems~\ref{T-genus3}
and~\ref{T-genus4} we began by searching for pointless curves
of genus $3$ and $4$ over various small finite fields.
In this section we explain the
heuristic we used to find families of curves
in which pointless curves might be abundant.
We begin with a lemma
from the theory of function fields over finite fields.

\begin{lem}
\label{L-Chebotarev}
Let $L/K$ be a degree-$d$ extension of function fields over 
a finite field~$k$, let $M$ be the Galois closure of $L/K$,
let $G = \Gal(M/K)$, and let $H = \Gal(M/L)$.
Let $S$ be the set of places $\frakp$ of $K$ that are
unramified in $L/K$ and for which there
is at least one place $\frakq$ of $L$, lying over $\frakp$,
with the same residue field as~$\frakp$.  Then the set $S$ 
has a Dirichlet density in the set of all places of $K$
unramified in $L/K$, and this density is
$$\delta:=\frac{\#\cup_{\tau\in G} H^\tau }{\# G}.$$
We have $\delta\ge 1/d$, with equality
precisely when $L$ is a Galois extension of~$K$.
Furthermore, we have $\delta\le 1 - (d-1)/\#G$.
\end{lem}

\begin{proof}
An easy exercise in the class field theory of function fields
({\it cf.}~\cite[proof of Lemma~2]{vdHeiden})
shows that the set $S$ is precisely the set of places $\frakp$
whose Artin symbol $(\frakp, L/K)$ lies in 
the union of the conjugates of $H$ in $G$.  The density statement
then follows from the Chebotarev density theorem.  

Since $H$ is an index-$d$ subgroup of $G$, we have
$$\frac{\#\cup_{\tau\in G} H^\tau }{\# G} \ge 
   \frac{\# H}{\# G} = \frac{1}{d}.$$
If $L/K$ is Galois then $H$ is trivial and the first relation
in the displayed equation above is an equality.
If $L/K$ is not Galois then $H$ is a nontrivial non-normal subgroup
of $G$, so the first relation above is an inequality.

To prove the upper bound on $\delta$, we note that two
conjugates $H^\sigma$ and $H^\tau$  of $H$ are identical when 
$\sigma$ and $\tau$ lie in the same coset of $H$ in $G$,
so when we form the union of the conjugates of $H$ we need
only let $\tau$ range over a set of coset representatives
of the $d$ cosets of $H$ in~$G$.
Furthermore, the identity element lies in every conjugate of~$H$,
so the union of the conjugates of $H$ contains at
most $d\cdot\#H - (d-1)$ elements. The upper bound follows.
\end{proof}

Note that the density mentioned in Lemma~\ref{L-Chebotarev} is 
a Dirichlet density.  If the constant field of $K$ is
algebraically closed in the Galois closure of $L/K$, then
the set $S$ also has a natural density
(see~\cite{MurtyScherk}).  In particular, the set $S$
has a natural density when $L/K$ is a Galois extension and 
$L$ and $K$ have the same constant field.

Lemma~\ref{L-Chebotarev} leads us to our main heuristic:

\begin{heur}
\label{H-automorphisms}
Let $C\to D$ be a degree-$d$ cover of curves over~$\F_q$,
let $L/K$ be the corresponding extension of function fields,
and let $\delta$ be the density from Lemma~{\rm\ref{L-Chebotarev}}.
If the constant field of the Galois closure of $L/K$ is equal to~$\F_q$,
then $C$ will be pointless with probability
$(1-\delta)^{\#D(\F_q)}$.
In particular, if $C\to D$ is a Galois cover,
then $C$ will be pointless with probability $(1-1/d)^{\#D(\F_q)}$.
\end{heur}

\begin{proof}[Justification]
Lemma~\ref{L-Chebotarev} makes it reasonable
to expect that with probability $1 - \delta$, 
a given rational point of $D$ will have no rational points of $C$
lying over it.  
Our heuristic follows if we assume that all 
of the points of $D$ behave independently.
\end{proof}

Consider what this heuristic tells us about hyperelliptic
curves.  Since a hyperelliptic curve is a double cover of a genus-$0$
curve,
we expect that a hyperelliptic curve over $\F_q$ 
will be pointless with probability $(1/2)^{q+1}$.  However,
if the hyperelliptic curve has more automorphisms than just the
hyperelliptic involution, it will be more likely to be 
pointless.  For instance, suppose $C$ is a hyperelliptic
curve whose automorphism group has order~$4$.
This automorphism group will give us a Galois cover $C\to\BP^1$ of
degree~$4$.  Then our heuristic suggests that $C$
will be pointless with probability $(3/4)^{q+1}$.

On the other hand, consider a generic smooth plane quartic $C$ over~$\F_q$.
A generic quartic has a $1$-parameter family
of non-Galois maps of degree $3$ to $\BP^1$.
For any one of these maps, the Galois group of its Galois closure is
the symmetric group on $3$ elements.  In this case, the density $\delta$
from Lemma~\ref{L-Chebotarev} is $2/3$, so we expect (modulo the condition
on constant fields mentioned in the heuristic) that a typical plane quartic 
will be pointless with probability $(1/3)^{q+1}$.
But if the quartic $C$ has an automorphism group
of order $4$, and if the quotient of $C$ by this automorphism group
is $\BP^1$, then we expect $C$ to be pointless with
probability $(3/4)^{q+1}$.

This heuristic suggested two things to us.  First, to find pointless
curves it is helpful to look for curves with larger-than-usual
automorphism groups.  We decided to focus on curves whose
automorphism groups contain the Klein $4$-group, because it is
easy to write down curves with this automorphism group and yet the
group is large enough to give us a good chance of finding
pointless curves.  
Second, the heuristic suggested that we look at curves $C$ that are
double covers of curves $D$ that are double covers of~$\BP^1$.
The Galois group of the resulting degree-$4$ cover $C\to\BP^1$ will typically
be the dihedral group of order~$8$, and the heuristic predicts that 
$C$ will be pointless with probability $(5/8)^{q+1}$. 
For a fixed $D$, if we consider the family of double covers $C\to D$
with $C$ of genus~$3$ or~$4$, our heuristic predicts that
$C$ will be pointless with probability $(1/2)^{\#D(\F_q)}$.
If $\#D(\F_q)$ is small enough, this probability can be reasonably high.  

The curves that we found by following our heuristic
are listed in Sections~\ref{S-examples3} and~\ref{S-examples4}.

\section{Proofs of the theorems}
\label{S-proof}

In this section we prove the theorems stated in the introduction.
Clearly Theorem~\ref{T-genus3} follows from
Theorem~\ref{T-genus3specific}, so we will only prove 
Theorems~\ref{T-genus4} and~\ref{T-genus3specific}.

\begin{proof}[Proof of Theorem~{\rm\ref{T-genus3specific}}]
The Weil bound says that a curve of genus $3$ over $\F_q$ has
at least $q + 1 - 6\sqrt{q}$ points, and it follows immediately
that if there is a pointless genus-$3$ curve over $\F_q$ then
$q < 33$.  In Section~\ref{S-examples3} we give examples
of pointless genus-$3$ hyperelliptic
curves over $\F_q$ for $q\le 25$ and examples
of pointless smooth plane quartics for $q\le 23$, for
$q = 29$, and for $q = 31$.
To complete the proof, 
we need only prove the following statements:
\begin{enumerate}
\item\label{31} There are no pointless genus-$3$ curves over $\F_{31}$.
\item\label{27} There are no pointless genus-$3$ curves over $\F_{27}$.
\item\label{25} There are no pointless smooth plane quartics over $\F_{25}$.
\item\label{32} There are no pointless genus-$3$ hyperelliptic curves
over~$\F_{32}$.
\item\label{29} There are no pointless genus-$3$ hyperelliptic curves
over~$\F_{29}$.
\end{enumerate}

\subsubsection*{Statement~\ref{31}}
Theorem~1 of~\cite{LauterSerre:CM} shows that every genus-$3$ curve
over $\F_{31}$ has at least $2$ rational points, and
statement~\ref{31} follows.

\subsubsection*{Statement~\ref{27}}
To prove statement~\ref{27}, we begin by
running the Magma program {\tt CheckQGN} described
in~\cite{HoweLauter}.  The output of {\tt CheckQGN(27,3,0)}
shows that if $C$ is a pointless genus-$3$
curve over $\F_{27}$ then the real Weil polynomial of $C$ 
(see~\cite{HoweLauter}) must be $(x-10)^2 (x-8)$.
(To reach this conclusion without relying on the computer,
one can adapt the reasoning on `defect 2' found 
in~\cite[\S~2]{LauterSerre:JAG}.)
Applying Proposition~13 of~\cite{HoweLauter}, we find that
$C$ must be a double cover of an elliptic curve over $\F_{27}$ with
exactly $20$ rational points.

Up to Galois conjugacy, there are two elliptic curves
over $\F_{27}$ with exactly $20$ rational points; one is given by
$y^2 = x^3 + 2x^2 + 1$ and the other by 
$y^2 = x^3 + 2x^2 + a$, where $a^3 - a + 1 = 0$.
By using the argument given in the analogous situation
in~\cite[\S~6.1]{HoweLauter}, we see that
every genus-$3$ double cover of one of these two $E$'s can be
obtained by adjoining to the function field of $E$ an element $z$
that satisfies $z^2 = f$, where $f$ is a function on $E$ of degree
at most $6$ that is regular outside $\infty$, that has four
zeros or poles of odd order, and that has a double zero at a 
point $Q$ of $E$ that is rational over~$\F_{27}$.
In fact, it suffices to consider $Q$'s that represent the
classes of $E(\F_{27}) / 2 E(\F_{27})$.  The first $E$ 
given above has four such classes and the second has two.
We can also demand that the representative points $Q$ not be
$2$-torsion points.

The divisor of the function $f$ is
$$P_1 + P_2 + P_3 + P_4 + 2Q - 6\infty$$
for some geometric points $P_1,\ldots,P_4$.  We are
assuming that the double cover $C$ has no rational points, 
so none of the $P_i$ can be rational over $\F_{27}$.  In particular,
none of the $P_i$ is equal to the infinite point.  Since $Q$ is
also not the infinite point (because we chose it not to be
a $2$-torsion point), we see that the degree of $f$ is exactly~$6$.

It is easy to have Magma enumerate, for each of the six $(E,Q)$
pairs, all of the degree-$6$ functions $f$ on $E$ that have double
zeros at~$Q$.  For each such $f$ we can check to see whether there 
is a rational point $P$ on $E$ such that $f(P)$ is a nonzero
square; if there is such a point, then the double $D$ cover of $E$
given by $z^2 = f$ would have a rational point.
For those functions $f$ for which such a $P$ does not exist,
we can check to see whether the divisor of $f$ has the right form.
If the divisor of $f$ does have the right form, we can compute
whether the curve $D$ has a rational point lying over $Q$ or 
over~$\infty$.  

We wrote Magma routines to perform these calculations; they
are available on the web at the URL mentioned in the
acknowledgments.  As it happens, no $(E,Q)$ pair gives
rise to a function $f$ that passes the first two tests
described in the preceding paragraph, so we never had to
perform the third test. 

Our conclusion is that there are no pointless genus-$3$ curves
over $\F_{27}$, which completes the proof of statement~\ref{27}.

\subsubsection*{Statement~\ref{25}}
To prove statement~\ref{25} we start by running {\tt CheckQGN(25,3,0)}.
We find that the real Weil polynomial of a pointless genus-$3$ 
curve over $\F_{25}$ is either 
$f_1 := (x - 10)^2(x - 6)$ or $f_2:=(x - 10)(x^2 - 16 x + 62)$
or $f_3:=(x - 10)(x - 9)(x - 7)$ or $f_4:=(x - 10)(x - 8)^2$.
(This list can also be obtained by using Table~4 and Theorem~1(a)
of~\cite{HoweLauter}.)

We begin by considering the real Weil polynomial $f_1 = (x-10)^2 (x-6)$.
Suppose $C$ is a genus-$3$ curve over $\F_{25}$ with
real Weil polynomial equal to~$f_1$.
Arguing as in the proof of~\cite[Cor.~12]{HoweLauter}, we find
that there is an exact sequence
$$0 \to \Delta \to A\times E \to \Jac C \to 0,$$
where $A$ is an abelian surface with real Weil polynomial
$(x-10)^2$, where $E$ is an elliptic curve with real 
Weil polynomial $x-6$,
where $\Delta$ is a self-dual finite group scheme that
is killed by $4$, and
where the projections from $A\times E$ to $A$ and to $E$
give monomorphisms $\Delta\hookrightarrow A$ and
$\Delta\hookrightarrow E$.
Furthermore, there are polarizations $\lambda_A$ and $\lambda_E$
on $A$ and $E$ whose 
kernels are the images of $\Delta$ under these monomorphisms,
and the polarization on $\Jac C$ induced by the product
polarization $\lambda_A\times\lambda_E$ is the canonical
polarization on $\Jac C$.

Since $\Delta$ is isomorphic to the kernel of $\lambda_E$ and
since $\Delta$ is killed by $4$, we see that if $\Delta$ is
not trivial then it is isomorphic to either $E[2]$ or $E[4]$.
If $\Delta$ were trivial then $\Jac C$ would be equal to 
$A\times E$ and the canonical polarization on $\Jac C$ would be
a product polarization, and this is not possible.  Therefore
$\Delta$ is isomorphic either to $E[2]$ or $E[4]$.  Since
the Frobenius endomorphism of $A$ is equal to the multiplication-by-$5$
map on $A$, the group of geometric $4$-torsion points on $A$ is a 
trivial Galois module.  But $E[4]$ is not a trivial Galois 
module, so we see that $\Delta$ must be isomorphic to $E[2]$.
Arguing as in the proof of~\cite[Prop.~13]{HoweLauter},
we find that there must be a degree-$2$ map from $C$ to~$E$.

Thus, to find the genus-$3$ curves over $\F_{25}$ whose real Weil
polynomials are equal to  $(x-10)^2(x-6)$, we need only look at the genus-$3$
curves that are double covers of elliptic curves over $\F_{25}$
with $20$ points and with three rational points of order~$2$.
There are two such elliptic curves, and, as in the proof of 
statement~\ref{27}, we can use Magma to enumerate
their genus-$3$ double covers with no points.
(Our Magma program is available at the URL mentioned in the
acknowledgments.)
We find that there is exactly one such double cover: if $a$ is an element
of $\F_{25}$ with $a^2 - a + 2 = 0$, then the double cover $C$ of
the elliptic curve $y^2 = x^3 + 2x$ given by setting $z^2 = a(x^2-2)$
has no points.

The curve $C$ is clearly hyperelliptic, because it is a double cover
of the genus-$0$ curve $z^2 = a(x^2-2)$.  By parametrizing this
genus-$0$ curve and manipulating the resulting equation for $C$, 
we find that $C$ is isomorphic to the curve $y^2 = a(x^8 + 1)$,
which is the example presented below in Section~\ref{S-examples3}.

Next we show that there are no pointless genus-$3$
curves over $\F_{25}$ with real Weil polynomial equal to 
$f_2$ or $f_3$ or~$f_4$.

Suppose $C$ is a pointless genus-$3$ curve over $\F_{25}$
whose real Weil polynomial is $f_2$ or $f_3$ or~$f_4$.
By applying Proposition~13 of~\cite{HoweLauter}, we find
that $C$ must be a double cover of an elliptic curve over $\F_{25}$
having either $16$ or $17$ points.  There is one elliptic curve over
$\F_{25}$ of each of these orders.  As we did above and in 
the proof of statement~\ref{27}, we can easily have
Magma enumerate the genus-$3$ double covers of these elliptic
curves.  
The only complication is that for the curve with $16$ points,
we cannot assume that the auxiliary point $Q$ mentioned in
the proof of statement~\ref{27} is not a $2$-torsion point.

The Magma program we used to enumerate these double covers
can be found at the web site mentioned in the acknowledgments.
Using this program, we found that the curve with $17$ points has no pointless
genus-$3$ double covers.  On the other hand, we found two functions $f$ on the
curve $E$ with $16$ points such that the double cover of $E$ defined by
$z^2 = f$ is a pointless genus-$3$ curve. But when we
computed an upper bound for the number of points on these
curves over $\F_{625}$, we found that both of the curves have
at most $540$ points over $\F_{625}$.  This upper bound is
not consistent with any of the three real Weil polynomials we
are considering. (In fact, one can show by direct computation that
the two curves are isomorphic to the curve  $y^2 = a(x^8 + 1)$ that
we found earlier, whose real Weil polynomial is~$f_1$.)
Thus, there are no pointless genus-$3$ curves
over $\F_{25}$ with real Weil polynomial equal to $f_2$ or $f_3$ 
or~$f_4$.

This proves statement~\ref{25}.

\subsubsection*{Statement~\ref{32}}
Suppose that $C$ is a pointless genus-$3$ curve over~$\F_{32}$.
If $C$ were hyperelliptic, then its quadratic twist would 
be a genus-$3$ curve over~$\F_{32}$ with $66$ rational points.
But~\cite[Thm.~1]{LauterSerre:JAG} shows that no such curve
exists.

We give a second proof of statement~\ref{32} as well, which
provides us with a little extra information and foreshadows
some of our later arguments.  This same proof is given 
in~\cite[\S~3.3]{Elkies} and attributed to Serre.

Suppose that $C$ is a pointless genus-$3$ curve over~$\F_{32}$.
Then $C$ meets the Weil-Serre lower bound, and (as Serre shows
in~\cite{Serre:notes}) its Jacobian
is therefore isogenous to the cube of an elliptic curve $E$
over $\F_{32}$ whose trace of Frobenius is~$11$.  Note that the
endomorphism ring of this elliptic curve is
the quadratic order $\calO$ of discriminant  $11^2 - 4\cdot32 = -7$.
The polarizations of abelian varieties isogenous to a power
of a single elliptic curve whose endomorphism ring is a
maximal order can be understood in terms of Hermitian modules
(see the appendix to~\cite{LauterSerre:CM}).   
Since the endomorphism ring $\calO$  
is a maximal order and a PID, there is exactly one 
abelian variety in the isogeny class of $E^3$, namely $E^3$
itself.  Furthermore, the theory of Hermitian modules
shows that the principal polarizations of $E^3$ correspond 
to the isomorphism classes of unimodular Hermitian forms on the 
$\calO$-module $\calO^3$.  Hoffmann~\cite{Hoffmann} shows that
there is only one isomorphism class of
indecomposable unimodular Hermitian forms on $\calO^3$,
so there is at most one Jacobian in the isogeny class of~$E^3$,
and hence
at most one genus-$3$ curve over $\F_{32}$ with no points.
The example we give in Section~\ref{S-examples3} is a plane quartic,
so there are no pointless genus-$3$ hyperelliptic curves over~$\F_{32}$.
This proves statement~\ref{32}.

\subsubsection*{Statement~\ref{29}}
We wrote a Magma program to find (by enumeration) all pointless
genus-$3$ hyperelliptic curves over an arbitrary finite field $\F_q$
of odd characteristic with $q>7$.  We applied our program to the
field~$\F_{29}$, and we found no curves.
Our Magma program is available at the URL mentioned
in the acknowledgments.
\end{proof}

Note that in the course of proving Theorem~\ref{T-genus3specific}
we showed that the pointless genus-$3$ curves over $\F_{25}$ and
$\F_{32}$ exhibited in Section~\ref{S-examples3} 
are the only such curves over their respective fields.
Also, our program to enumerate pointless
genus-$3$ hyperelliptic curves shows that there is only one
pointless genus-$3$ hyperelliptic curve over~$\F_{23}$.

\begin{proof}[Proof of Theorem~{\rm\ref{T-genus4}}]
It follows from Serre's refinement of the Weil
bound~\cite[Th\'e\-or\`eme~1]{Serre:CRAS}
that if a curve of genus $4$ over $\F_q$ has no rational points, 
then~$q\le 59$.
In Section~\ref{S-examples4} we give examples
of pointless genus-$3$ curves over $\F_q$ for all $q$ with~$q\le 49$,
so to prove the theorem we must show that there are no pointless
genus-$4$ curves over $\F_{53}$ or~$\F_{59}$.

Combining the output of {\tt CheckQGN(53,4,0)} with Theorem~1(b)
of~\cite{HoweLauter}, we find that a pointless genus-$4$ curve 
over~$\F_{53}$ must be a double cover of an elliptic curve $E$
over~$\F_{53}$ with exactly~$42$ points.  (Again, the information
obtained by running {\tt CheckQGN} can also be
obtained without recourse to the computer by modifying
the `defect $2$' arguments in~\cite[\S~2]{LauterSerre:JAG}.)

There are four elliptic curves $E$ over $\F_{53}$ with
exactly $42$ points.  Following the arguments
of~\cite[\S~6.1]{HoweLauter}, we find that every
genus-$4$ double cover of such an $E$ can be obtained by adjoining 
to the function field of $E$ a root of an equation
$z^2 = f$, where $f$ is a function on $E$ whose divisor
is of the form
$$P_1 + \cdots + P_6 + 2Q - 8\infty,$$
where $Q$ is a rational point of $E$ that is not killed by $2$,
and where it suffices to consider $Q$ that cover the
residue classes of $E(\F_{53})$ modulo $3 E(\F_{53})$.
As in the preceding proof, we wrote Magma programs to
enumerate the genus-$4$ double covers of the four possible $E$'s
and to check to see whether all of these covers had rational
points.  Our programs, available at the URL mentioned in the
acknowledgments, showed that every genus-$4$ double
cover of these $E$'s has a rational point.  Thus there are
no pointless genus-$4$ curves over $\F_{53}$.

Next we show that there are no pointless curves of genus $4$ 
over $\F_{59}$.  If $C$ were such a curve, then $C$ would meet
the Weil-Serre lower bound, and therefore  the Jacobian of $C$ 
would be isogenous
to the fourth power of an elliptic curve $E$ over $\F_{59}$
with $45$ points.  Note that there is exactly one such $E$,
and its endomorphism ring $\calO$ is the quadratic order
of discriminant~$-11$.  As in the proof of statement~\ref{32}
of the proof of Theorem~\ref{T-genus3specific},
we see that there is
only one abelian variety in the isogeny class of $E^4$,
and principal polarizations of $E^4$ correspond 
to the isomorphism classes of unimodular Hermitian forms on the 
$\calO$-module~$\calO^4$.  
Schiemann~\cite{Schiemann}
states that there are six isomorphism classes of unimodular Hermitian
forms on the module~$\calO^4$.  We were unable to find a listing of
these isomorphism classes at the URL mentioned in~\cite{Schiemann},
but we did find them by following links from the URL

\smallskip
\noindent
{\tt http://www.math.uni-sb.de/\~{}ag-schulze/Hermitian-lattices/}
\smallskip

\noindent
We have put a copy of the page listing these six forms on the
web site mentioned in the acknowledgments.

Three of the isomorphism classes of unimodular Hermitian forms on $\calO^4$
are decomposable, and so do not come from the Jacobian of a 
curve.  The three indecomposable Hermitian forms can each be
written as a matrix with an upper left entry of~$2$.
Arguing as in the proof of~\cite[Prop.~13]{HoweLauter},
we find that our curve $C$ must be a double cover of the curve~$E$.

We are again in familiar territory.  As above, it is an easy matter
to write a Magma program to enumerate the genus-$4$ double covers
of the given elliptic curve $E$ and to check that they all have
a rational point.  (Our Magma programs are available at the URL
mentioned in the acknowledgments.)  Our computation showed that
there are no pointless curves of genus $4$ over~$\F_{59}$.
\end{proof}

\section{Examples of pointless curves of genus $3$}
\label{S-examples3}

In this section we give examples of pointless curves of
genus $3$ over the fields where such curves exist.  
We only consider curves whose automorphism groups contain
the Klein $4$-group~$V$.
We begin with the hyperelliptic curves.

Suppose $C$ is a genus-$3$ hyperelliptic curve over $\F_q$
whose automorphism group contains a copy of~$V$, and assume that the
hyperelliptic involution is contained in~$V$.
Then $V$ modulo the hyperelliptic involution acts on $C$ modulo
the hyperelliptic involution, and gives us an involution on~$\BP^1$.
By changing coordinates on $\BP^1$, we may assume that the
involution on $\BP^1$ is of the form $x\mapsto n/x$ for some
$n\in \F_q^*$.  (When $q$ is odd we need consider only two values 
of~$n$, one a square and one a nonsquare.  When $q$ is even 
we may take $n=1$.)

It follows that when $q$ is odd the curve $C$ can be defined
either by an equation of the form $y^2 = f(x + n/x)$, where $f$ is a
separable quartic polynomial coprime to $x^2 - 4n$,
or by an equation of the form $y^2 = x f(x + n/x)$, where $f$ is a 
separable cubic polynomial coprime to $x^2 - 4n$.  However, the latter
possibility cannot occur if $C$ is to be pointless.
When $q$ is even, if we assume the curve if ordinary then 
it may be written in the form $y^2 + y = f(x + 1/x)$, where
$f$ is a rational function with $2$ simple poles, both nonzero.

We wrote a simple Magma program to search for pointless
hyperelliptic curves of this form.  We found such curves for
every $q$ in 
$$\{2, 3, 4, 5, 7, 8, 9, 11, 13, 16, 17, 19, 23, 25\}.$$
We give examples in Table~\ref{Tbl-examples3}.

\begin{table}
\renewcommand{\arraystretch}{1.25}
\begin{center}
\begin{tabular}{|r|l|}
\hline
$q$ & curve \\ \hline  \hline
$2$ &
$y^2 + y = (x^4 + x^2 + 1)/(x^4 + x^3 + x^2 + x + 1)$
\\ \hline
$3$ &
$y^2 = - x^8 + x^7 - x^6 - x^5 - x^3 - x^2 + x - 1$
\\ \hline
$4$ &
$y^2 + y = (ax^4 + ax^3 + a^2x^2 + ax + a)/(x^4 + ax^3 + x^2 + ax + 1)$\\
& where $a^2 + a + 1 = 0$ \\ \hline
$5$ &
$y^2 = 2x^8 + 3x^4 + 2$
\\ \hline
$7$ &
$y^2 = 3x^8 + 2x^6 + 3x^4 + 2x^2 + 3$
\\ \hline
$8$ &
$y^2 + y = (x^4 + a^6x^3 + a^3x^2 + a^6x + 1)/(x^4 + x^3 + x^2 + x + 1)$ 
\\ & where $a^3 + a + 1 = 0$\\ \hline
$9$ &
$y^2 = a(x^8 + 1)$
\\ & where $a^2 - a - 1= 0$ \\ \hline
$11$ &
$y^2 = 2x^8 + 4x^6 - 2x^4 + 4x^2 + 2$
\\ \hline
$13$ &
$y^2 = 2x^8 + 3x^7 + 3x^6 + 4x^4 + 3x^2 + 3x + 2$
\\ \hline
$16$ &
$y^2 + y = (a^3x^4 + a^3x^3 + a^{14}x^2 + a^3x + a^3)/(x^4 + a^3x^3 + x^2 + a^3x + 1)$
\\ & where $a^4 + a + 1 = 0$\\ \hline
$17$ &
$y^2 = 3x^8 - 2x^5 + 4x^4 - 2x^3 + 3$
\\ \hline
$19$ &
$y^2 = 2x^8 - x^6 - 8x^4 - x^2 + 2$
\\ \hline
$23$ & 
$y^2 = 5x^8 + x^6 + 6x^5 + 7x^4 - 6x^3 + x^2 + 5$
\\ \hline
$25$ &
$y^2 = a(x^8 + 1)$
\\ & where $a^2 - a + 2 = 0$\\ \hline
\end{tabular}
\end{center} 
\vspace{1ex}
\caption{
Examples of pointless hyperelliptic curves of genus $3$ over $\F_q$
with automorphism group containing the Klein $4$-group.  For
$q\neq 23$, the automorphism $x \mapsto 1/x$ of $\BP^1$ lifts to give
an automorphism of the curve; for $q=23$, the automorphism
$x\mapsto -1/x$ lifts.}
\label{Tbl-examples3} 
\end{table}

Now we turn to the pointless smooth plane quartics.
We searched for pointless quartics of the form 
$$ax^4 + by^4 + cz^4 + dx^2y^2 + ex^2z^2 + fy^2z^2 = 0$$
over finite fields of odd characteristic, 
because the automorphism groups of such quartics clearly contain
the Klein group.  We found pointless quartics of this form
over $\F_q$ for $q$ in
$$\{5,7,9,11,13,17,19,23,29\}.$$  
We present sample curves in Table~\ref{Tbl-quartics}.

\begin{table}
\renewcommand{\arraystretch}{1.25}
\begin{center}
\begin{tabular}{|r|l|}
\hline
$q$ & curve \\ \hline  \hline
$5$ &
$x^4 + y^4 + z^4 = 0$
\\ \hline
$7$ &
$x^4 + y^4 + 2z^4 + 3x^2z^2 + 3y^2z^2 = 0$
\\ \hline
$9$ &
$x^4 - y^4 + a^2z^4 + x^2y^2 = 0$
\\ & where $a^2 - a - 1 = 0$ \\ \hline
$11$ &
$x^4 + y^4 + z^4 + x^2y^2 + x^2z^2 + y^2z^2 = 0$
\\ \hline
$13$ &
$x^4 + y^4 + 2z^4 = 0$
\\ \hline
$17$ &
$x^4 + y^4 + 2z^4 + x^2y^2 = 0$
\\ \hline
$19$ &
$x^4 + y^4 + z^4 + 7x^2y^2 - x^2z^2 - y^2z^2 = 0$
\\ \hline
$23$ & 
$x^4 + y^4 + z^4 + 10x^2y^2 - 3x^2z^2 - 3y^2z^2 = 0$
\\ \hline
$29$ &
$x^4 + y^4 + z^4 = 0$
\\ \hline
\end{tabular}
\end{center} 
\vspace{1ex}
\caption{
Examples of pointless smooth plane quartics over $\F_q$ (with $q$ odd)
with automorphism group containing the Klein $4$-group.}
\label{Tbl-quartics} 
\end{table}

Over $\F_3$ there are many pointless smooth plane quartics; for instance,
the curve $$ x^4 + xyz^2 + y^4 + y^3z - yz^3 + z^4 = 0 $$ has no points. 

We know from the proof of Theorem~\ref{T-genus3specific} that 
there is at most one pointless genus-$3$ curve over $\F_{32}$, 
and its Jacobian is isomorphic to the cube of an elliptic curve
whose endomorphism ring has discriminant~$-7$.  This suggests that
we should look at twists of the reduction of the Klein quartic,
and indeed we find that the curve 
$$(x^2 + x)^2 + (x^2 + x)(y^2 + y) + (y^2 + y)^2 + 1 = 0$$
has no points over $\F_{32}$.
(This fact is noted in~\cite[\S~3.3]{Elkies}.)
For the other fields of characteristic $2$, we 
find examples by modifying the example for~$\F_{32}$.
We list the results in 
Table~\ref{Tbl-quartics2}.

\begin{table}
\renewcommand{\arraystretch}{1.25}
\begin{center}
\begin{tabular}{|r|l|}
\hline
$q$ & curve \\ \hline  \hline
$2$ &
$(x^2 + xz)^2 + (x^2 + xz)(y^2 + yz) + (y^2 + yz)^2 + z^4 = 0$
\\ \hline
$4$ &
$(x^2 + xz)^2 + a(x^2 + xz)(y^2 + yz) + (y^2 + yz)^2 + a^2z^4 = 0$
\\ & where $a^2 + a + 1 = 0$ \\ \hline
$8$ &
$(x^2 + xz)^2 + (x^2 + xz)(y^2 + yz) + (y^2 + yz)^2 + a^3z^4 = 0$
\\ & where $a^3 + a + 1 = 0$ \\ \hline
$16$ &
$(x^2 + xz)^2 + a(x^2 + xz)(y^2 + yz) + (y^2 + yz)^2 + a^7z^4 = 0$
\\ & where $a^4 + a + 1 = 0$ \\ \hline
$32$ &
$(x^2 + xz)^2 + (x^2 + xz)(y^2 + yz) + (y^2 + yz)^2 + z^4 = 0$
\\ \hline
\end{tabular}
\end{center} 
\vspace{1ex}
\caption{
Examples of pointless smooth plane quartics over $\F_q$ (with $q$ even)
with automorphism group containing the Klein $4$-group.}
\label{Tbl-quartics2} 
\end{table}

We close this section by mentioning a related method of 
constructing pointless genus-$3$ curves. 
Suppose $C$ is a genus-$3$ curve over a field of characteristic
not~$2$, and suppose that $C$ has a pair of commuting involutions
(like the curves we considered in this section).
Then either $C$ is an unramified double cover of a genus-$2$
curve, or $C$ is a genus-$3$ curve of the type considered 
in~\cite[\S~4]{HoweLeprevostPoonen}, that is, a genus-$3$
curve obtained by `gluing' three elliptic curves together
along portions of their $2$-torsion.
This suggests a more direct method of constructing genus-$3$
curves with no points:  We can start with three elliptic curves
with few points, and try to glue them together using the
construction from~\cite[\S~4]{HoweLeprevostPoonen}.
This idea was used by the third author to construct genus-$3$
curves with many points~\cite{Top}.

\section{Examples of pointless curves of genus $4$}
\label{S-examples4}

We searched for pointless genus-$4$ curves by looking
at hyperelliptic curves whose automorphism group
contained the Klein $4$-group; however, we found that 
for $q>31$ no such curves exist.  Since we need to find pointless genus-$4$ curves
over $\F_q$ for every $q\le 49$, we moved on to a different family
of curves with commuting involutions.

Suppose $q$ is an odd prime power
and suppose $f$ and $g$ are separable cubic polynomials
in $\F_q[x]$ with no factor in common.   An easy ramification 
computation shows that then the curve defined by $y^2 = f$ and $z^2 = g$
has genus~$4$.  Clearly the automorphism group of this
curve contains a copy of the Klein $4$-group.  
It is easy to check whether a curve of this form
is pointless:  For every value of $x$ in $\F_q$, at least one
of $f(x)$ and $g(x)$ must be a nonsquare, and exactly one of $f$ and
$g$ should have a nonsquare as its coefficient of $x^3$.
We found pointless curves of this form over every $\F_q$ with 
$q$ odd and $q\le 49$.   Examples are given in Table~\ref{Tbl-examples4}.

\begin{table}
\renewcommand{\arraystretch}{1.25}
\begin{center}
\begin{tabular}{|r|ll|}
\hline
$q$ & curve & \\ \hline  \hline
$3$ & 
$ y^2 = x^3 - x - 1$ &
$ z^2 = - x^3 + x - 1$
\\ \hline
$5$ & 
$ y^2 = x^3 - x + 2$ &
$ z^2 = 2 x^3 - 2 x$ 
\\ \hline
$7$ &
$ y^2 = x^3 - 3$ &
$ z^2 = 3 x^3 - 1$
\\
\hline $9$ & 
$ y^2 = x^3 - x + 1 $ &
$ z^2 = a (x^3 - x - 1) $
\\
& where $a^2 - a - 1 = 0$  &
\\ \hline
$11$ & 
$ y^2 = x^3 - x - 3$ &
$ z^2 = 2 x^3 - 2 x - 5$
\\ \hline
$13$ & 
$ y^2 = x^3 + 1 $ &
$ z^2 = 2 x^3 - 5$
\\ \hline
$17$ & 
$ y^2 = x^3 + x $ &
$ z^2 = 3 x^3 - 8 x^2 - 3 x + 5 $
\\ \hline
$19$ & 
$ y^2 = x^3 + 2$ &
$ z^2 = 2 x^3 + 1$
\\ \hline
$23$ & 
$ y^2 = x^3 + x + 6 $ &
$ z^2 = 5 x^3 + 9 x^2 - 3 x + 10$
\\ \hline
$25$ & 
$ y^2 = x^3 + x + 1 $ &
$ z^2 = a(x^3 + x^2 + 2)$
\\ 
& where $a^2 - a + 2 = 0$ &
\\ \hline
$27$ & 
$ y^2 = x^3 - x + a^5$ &
$ z^2 = -x^3 + x + a^5$ 
\\
& where $a^3 - a + 1 = 0$ &
\\ \hline
$29$ & 
$ y^2 = x^3 + x $ &
$ z^2 = 2 x^3 + 12 x + 14$
\\ \hline
$31$ & 
$ y^2 = x^3 - 10$ &
$ z^2 = 3 x^3 + 9$
\\ \hline
$37$ & 
$ y^2 = x^3 + x + 4 $ &
$ z^2 = 2 x^3 - 17 x^2 + 5 x + 15 $
\\ \hline
$41$ & 
$ y^2 = x^3 + x + 17 $ &
$ z^2 = 3 x^3 - x^2 - 12 x - 16$
\\ \hline
$43$ & 
$ y^2 = x^3 - 9$ &
$ z^2 = 2 x^3 + 18$
\\ \hline
$47$ & 
$ y^2 = x^3 + 5 x - 12 $ &
$ z^2 = 5 x^3 + 2 x^2 + 19 x - 9$
\\ \hline
$49$ & 
$ y^2 = x^3 + 4 $ &
$ z^2 = a (x^3 + 2) $
\\
& where $a^2 - a + 3 = 0$ &
\\ \hline
\end{tabular}
\end{center} 
\vspace{1ex}
\caption{
Examples of pointless curves of genus $4$ over $\F_q$ (with $q$ odd)
with automorphism group containing the Klein $4$-group.}
\label{Tbl-examples4} 
\end{table}

We mention two points of interest about curves of this form.
First, if the $\F_q$-vector subspace of $\F_q[x]$ spanned by
the cubic polynomials $f$ and $g$ contains the constant polynomial~$1$,
then the
curve $C$ defined by the two equations $y^2=f$ and $z^2=g$ is trigonal:
If we have $af + bg = 1$, then $(x,y,z)\mapsto (y,z)$ 
defines a degree-$3$ map from
$C$ to the genus-$0$ curve $ay^2 + bz^2 = 1$.
Second, if $q\equiv1\bmod 3$ and if the coefficients of $x$ and $x^2$ in
$f$ and $g$ are zero, then the curve $C$ has even more automorphisms,
given by multiplying $x$ by a cube root of unity.  (Likewise, if
$q$ is a power of $3$ and if $f$ and $g$ are both of the form
$a(x^3 - x) + b$, then $x\mapsto x+1$ gives an automorphism of~$C$.)
When it was possible, we chose the examples in Table~\ref{Tbl-examples4}
to have these properties. 
In Table~\ref{Tbl-examples4trigonal} we provide
trigonal models for the curves in Table~\ref{Tbl-examples4} that
have them.

\begin{table}
\renewcommand{\arraystretch}{1.25}
\begin{center}
\begin{tabular}{|r|l|l|}
\hline
$q$ & curve & liftable involutions of $\BP^1$ \\ \hline  \hline
$3$ &
$v^3 - v = (u^4 + 1)/(u^2 + 1)^2$ &
$u\mapsto -u,\quad u\mapsto 1/u$
\\ \hline
$5$ &
$v^3 - v = -2(u^2 - 2)^2/(u^2 + 2)^2$ &
$u\mapsto -u,\quad u\mapsto 2/u$
\\ \hline
$7$ &
$v^3 = 2 u^6 + 2$ &
$u\mapsto -u,\quad u\mapsto 1/u$
\\ \hline
$9$ &
$v^3 - v = (u^4 + a^2)/(u^2 + a^5)^2$ &
$u\mapsto -u,\quad u\mapsto a/u$
\\
& where $a^2 - a - 1 = 0$  &
\\ \hline
$11$ &
$v^3 - v = (3 u^4 + 4 u^2 + 3)/(u^2 + 1)^2$ &
$u\mapsto -u,\quad u\mapsto 1/u$
\\ \hline
$13$ &
$v^3 = 4 u^6 + 6$ &
$u\mapsto -u,\quad u\mapsto 2/u$
\\ \hline
$19$ &
$v^3 = 2 u^6 + 2$ &
$u\mapsto -u,\quad u\mapsto 1/u$
\\ \hline
$27$ &
$v^3 - v = a^{18} (u^4 + 1)/(u^2 + 1)^2$ &
$u\mapsto -u,\quad u\mapsto 1/u$
\\
& where $a^3 - a + 1 = 0$ &
\\ \hline
$31$ &
$v^3 = 5 u^6 - 11 u^4 - 11 u^2 + 5$ &
$u\mapsto -u,\quad u\mapsto 1/u$
\\ \hline
$43$ &
$v^3 = 7 u^6 + 8 u^4 + 8 u^2 + 7$ &
$u\mapsto -u,\quad u\mapsto 1/u$
\\ \hline
$49$ &
$v^3 = 2 u^6 + a$ &
$u\mapsto -u,\quad u\mapsto a^3/u$
\\
& where $a^2 - a + 3 = 0$ &
\\ \hline
\end{tabular}
\end{center} 
\vspace{1ex}
\caption{
Trigonal forms for some of the curves in Table~\ref{Tbl-examples4}.
The third column gives two involutions of $\BP^1$ that 
lift to give commuting involutions of the curve.}
\label{Tbl-examples4trigonal} 
\end{table}

It remains for us to find examples of pointless genus-$4$ curves
over $\F_2, \F_4, \F_8, \F_{16},$ and $\F_{32}$.

Let $q$ be a power of $2$.
An easy argument shows that 
a genus-$4$ hyperelliptic curve over $\F_q$ provided with an action of
the Klein group must have a rational Weierstra\ss\ point, and
so will not be pointless.  Thus we decided simply to enumerate 
the genus-$4$ hyperelliptic curves (with no rational
Weierstra\ss\ points) over the remaining $\F_q$ and 
to check for pointless curves.  We found pointless hyperelliptic
curves over $\F_q$ for $q\in \{2,4,8,16\}$;
the examples we give in Table~\ref{Tbl-genus4char2} are
all twists over $\F_q$ of curves that can be defined over~$\F_2$.

\begin{table}
\renewcommand{\arraystretch}{1.25}
\begin{center}
\begin{tabular}{|r|l|}
\hline
$q$ & curve \\ \hline  \hline
$2$ &
$y^2 + y = t + (x^4 + x^3 + x^2 + x)/(x^5 + x^2 + 1)$
\\ \hline
$4$ &
$y^2 + y = t + (x^3 + 1)/(x^5 + x^2 + 1)$
\\ \hline
$8$ &
$y^2 + y = t + (x^4 + x^3 + x^2 + x)/(x^5 + x^2 + 1)$
\\ \hline
$16$ &
$y^2 + y = t + (x^3 + 1)/(x^5 + x^2 + 1)$
\\ \hline
\end{tabular}
\end{center} 
\vspace{1ex}
\caption{
Examples of pointless genus-$4$ hyperelliptic curves
over $\F_q$ (with $q$ even).
On each line, the symbol $t$ refers to an arbitrary element
of $\F_q$ whose trace to $\F_2$ is equal to~$1$.}
\label{Tbl-genus4char2} 
\end{table}

Our computer search also revealed that every
genus-$4$ hyperelliptic curve over $\F_{32}$
has at least one rational point.
So to find an example of a pointless genus-$4$ curve over $\F_{32}$,
we decided to look for genus-$4$ double covers of elliptic curves~$E$.
Our heuristic suggested that we might have good luck finding pointless
curves if $E$ had few points, but for the sake of
completeness we examined every $E$ over~$\F_{32}$.

We found that up to isomorphism and Galois conjugacy
there are exactly two pointless genus-$4$ curves
over $\F_{32}$ that are double covers of elliptic curves.
The first can be defined by the equations
\begin{align*}
y^2 + y   & = x + 1/x + 1 \\
z^2 + z   & = \frac{a^7 x^4 + a^{30} x^3 y  + a^{13} x^2 
                            + x + a^{23} x y + a^6}
              {x^3 + a^{15} x^2 + x + a^{28}}
\end{align*}
and the second by
\begin{align*}
y^2 + y   & = x + a^7/x \\
z^2 + z   & = \frac{a^4 x^4 + a^7 x^3 y + a^3 x^3  
                       + a^{23} x^2 y + a^{28} x^2 + a^{28} x y + a^{16}}
              {x^3 + a^{25} x^2 + a^{22} x + a^{25}},
\end{align*}
where $a^5 + a^2 + 1 = 0.$



\begin{thebibliography}{99}    

\bibitem{AuerTop}
{\sc Roland Auer and Jaap Top}:
Some genus $3$ curves with many points,
pp.~163--171 in 
{\it Algorithmic Number Theory\/} 
(Claus Fieker and David R. Kohel, eds.),
Lecture Notes in Comp. Sci. {\bf 2369},
Springer-Verlag, Berlin, 2002. 

\bibitem{magma}
{\sc Wieb Bosma, John Cannon, and Catherine Playoust}:
The Magma algebra system. I. The user language.
{\it J. Symbolic Comput.} {\bf 24}  (1997) 235--265. 

\bibitem{Elkies}
{\sc Noam D. Elkies}:
The Klein quartic in number theory,
pp. 51--101 in:
{\it The eightfold way\/} (Silvio Levy, ed.),
Math. Sci. Res. Inst. Publ. {\bf 35},
Cambridge Univ. Press, Cambridge, 1999. 


\bibitem{GeerVlugt}
{\sc Gerard van der Geer and Marcel van der Vlugt}:
Tables of curves with many points,
{\it  Math. Comp.} {\bf 69}  (2000) 797--810. 
Updated versions available at {\tt http://www.science.uva.nl/\~{}geer/}.

\bibitem{Hasse}
{\sc H. Hasse}:
Zur Theorie der abstrakten elliptischen Funktionk\"orper. I, II, III,
{\it J. Reine Angew. Math.} {\bf 175} (1936) 55--62, 69--88, 193--208.

\bibitem{vdHeiden}
{\sc Gert-Jan van der Heiden}:
Local-global problem for Drinfeld modules,
{\it J. Number Theory} {\bf 104} (2004) 193--209.

\bibitem{Hoffmann}
{\sc Detlev W. Hoffmann}:
On positive definite Hermitian forms,
{\it Manuscripta Math.} {\bf 71} (1991) 399--429.

\bibitem{HoweLauter}
{\sc Everett W. Howe and Kristin E. Lauter}:
Improved upper bounds for the number of points on curves
over finite fields,
{\it Ann. Inst. Fourier {\rm(}Grenoble\/{\rm)}} {\bf 53} (2003) 1677--1737.
{\tt arXiv:math.NT/0207101}.

\bibitem{HoweLeprevostPoonen}
{\sc Everett W. Howe, Franck Lepr\'evost, and Bjorn Poonen}: 
Large torsion subgroups of split Jacobians of curves
of genus two or three,
{\it Forum Math.} {\bf 12} (2000) 315--364.

\bibitem{LauterSerre:JAG} 
{\sc Kristin Lauter with an appendix by J.-P. Serre}:
Geometric methods for improving the upper bounds on the
number of rational points on algebraic curves over finite fields, 
{\it J. Algebraic Geom.} {\bf 10} (2001) 19--36.
{\tt arXiv:math.AG/0104247}.

\bibitem{LauterSerre:CM}
{\sc Kristin Lauter with an appendix by Jean-Pierre Serre}:
The maximum or minimum number of rational points on 
genus three curves over finite fields, 
{\it Compositio Math.} {\bf 134} (2002) 87--111.
{\tt arXiv:math.AG/0104086}.

\bibitem{LeepYeomans}
{\sc David B. Leep and Charles C. Yeomans}:
Quintic forms over $p$-adic fields,
{\it J. Number Theory} {\bf 57} (1996) 231--241. 

\bibitem{MaisnerNart}
{\sc Daniel Maisner and Enric Nart with an appendix by Everett W. Howe}: 
Abelian surfaces over finite fields as Jacobians, 
{\it Experiment. Math.} {\bf 11} (2002) 321Ð337.

\bibitem{MurtyScherk}
{\sc Vijaya Kumar Murty and John Scherk}:
Effective versions of the Chebotarev density theorem for function fields,
{\it C. R. Acad. Sci. Paris S\'er. I Math.} {\bf 319} (1994) 523--528.

\bibitem{Schiemann}
{\sc Alexander Schiemann}:
Classification of Hermitian forms with the neighbour method,
{\it J. Symbolic Computation} {\bf 26} (1998) 487--508.

\bibitem{Serre:CG}
{\sc Jean-Pierre Serre}:
{\it Cohomologie Galoisienne 
(cinqui\`eme \'edition, r\'evis\'ee et compl\'et\'ee)},
Lecture Notes in Math.~{\bf 5},
Springer-Verlag, Berlin, 1994.

\bibitem{Serre:CRAS}
{\sc Jean-Pierre Serre}:
Sur le nombre des points rationnels d'une courbe 
alg\'ebrique sur un corps fini, 
{\it C. R. Acad. Sci. Paris S\'er. I Math.} {\bf 296} 
(1983) 397--402; = \OE{}uvres [128].

\bibitem{Serre:notes}
{\sc Jean-Pierre Serre}:
{\it Rational points on curves over finite fields},
unpublished notes by Fernando Q. Gouv\^ea of lectures 
at Harvard University, 1985.

\bibitem{Stark}
{\sc H. M. Stark}:
On the Riemann hypothesis in hyperelliptic function fields,
pp. 285--302 in:
{\it Analytic number theory\/} (Harold G. Diamond, ed.),
Proc. Sympos. Pure Math. {\bf 24},
American Mathematical Society, Providence, R.I. 1973.

\bibitem{Top}
{\sc Jaap Top}:
Curves of genus $3$ over small finite fields,
{\it Indag. Math. (N.\,S.)} {\bf 14} (2003) 275--283.

\end{thebibliography}
\end{document}